\documentclass[12pt]{amsart}
\usepackage{amssymb}
\usepackage{amsmath}
\usepackage{amsfonts}
\usepackage{mathrsfs}
\usepackage{graphicx}
\usepackage{color}
\usepackage[onehalfspacing]{setspace}
\usepackage{caption}
\usepackage{natbib}
\usepackage{enumerate}
\usepackage[utf8]{inputenc}
\usepackage{charter}
\usepackage[colorlinks=true,citecolor=blue,urlcolor=blue,pdfpagemode=UseNone,pdfstartview=FitH]{hyperref}

\makeatletter
\def\section{\@startsection{section}{1}
	\z@{1.0\linespacing\@plus\linespacing}{.8\linespacing}{\Large}}

\def\subsection{\@startsection{subsection}{2}
	\z@{.8\linespacing\@plus.7\linespacing}{.7\linespacing}{\large}}

\def\subsubsection{\@startsection{subsubsection}{3}
	\z@{.5\linespacing\@plus.7\linespacing}{-.5em}{\normalfont\bfseries}}
\makeatother

\newtheorem{theorem}{Theorem}[section]
\newtheorem{proposition}{Proposition}[section]
\newtheorem{lemma}{Lemma}[section]
\newtheorem{corollary}{Corollary}[section]

\theoremstyle{definition}
\newtheorem{definition}{Definition}[section]

\theoremstyle{definition}

\theoremstyle{definition}

	\title{}
	
\begin{document}
		\vspace*{5ex minus 1ex}
		\begin{center}
			\LARGE \textsc{A Uniform-in-$P$ Edgeworth Expansion under Weak Cram\'{e}r Conditions}
			\bigskip
		\end{center}
		
		\date{%
			\today%
		}
		
		\vspace*{3ex minus 1ex}
		\begin{center}
			Kyungchul Song\\
			\textit{Vancouver School of Economics, University of British Columbia}\\
			\bigskip
			\bigskip
		\end{center}
		
		\thanks{I thank three anonymous referees for valuable comments in a previous version of this paper. All errors are mine. I acknowledge that this research was supported by Social Sciences and Humanities Research Council of Canada. Corresponding address: Kyungchul Song, Vancouver School of Economics, University of British Columbia, Vancouver, BC, Canada. Email address: kysong@mail.ubc.ca.}
		\address{Vancouver School of Economics, University of British Columbia, 6000 Iona Drive, Vancouver, BC, V6T 1L4, Canada}
		\email{kysong@mail.ubc.ca}
		
		\fontsize{12}{14} \selectfont
		
\begin{abstract}
 This paper provides a finite sample bound for the error term in the Edgeworth expansion for a sum of independent, potentially discrete, nonlattice random vectors, using a uniform-in-$P$ version of the weaker Cram\'{e}r condition in \cite{Angst/Poly:17:EJP}. This finite sample bound can be used to derive an Edgeworth expansion that is uniform over the distributions of the random vectors. Using this result, we derive a uniform-in-$P$ higher order expansion of resampling-based distributions.
\medskip

{\noindent \textsc{Key words.} Edgeworth Expansion; Normal Approximation; Bootstrap; Resampling; Weak Cram\'{e}r Condition}
\medskip

{\noindent \textsc{AMS 1991 Classification: 60E05, 62G20, 62E20}}
\end{abstract}

\maketitle

\section{Introduction}
Suppose that $\{X_{i,n}\}_{i=1}^n$ is a triangular array of independent random vectors taking values in $\mathbf{R}^d$ and having mean zero. Let $\mathscr{P}_n$ be the collection of the joint distributions for $(X_{i,n})_{i=1}^n$. (One can view each joint distribution $P$ in $\mathscr{P}_n$ as being induced by a probability measure on a measurable space on which $(X_{i,n})_{i=1}^n$ live.) Define
\begin{eqnarray*}
	V_{n,P} \equiv \frac{1}{n}\sum_{i=1}^n \text{Var}_P(X_{i,n}),
\end{eqnarray*}
where $\text{Var}_P(X_{i,n}) \equiv \mathbf{E}_P(X_{i,n}X_{i,n}')$ and $\mathbf{E}_P$ denotes the expectation with respect to $P$. We are interested in the Edgeworth expansion of the distribution of
\begin{eqnarray*}
	\frac{1}{\sqrt{n}}\sum_{i=1}^n V_{n,P}^{-1/2}X_{i,n},
\end{eqnarray*}
which is uniform over $P \in \mathscr{P}_n$.

The Edgeworth expansion has long received attention in the literature. See \cite{Bhattacharya/RangaRao:10:NormAppprox} for a formal review of the results. The validity of the classical Edgeworth expansion is obtained under the Cram\'{e}r condition which says that the average of the characteristic functions over the sample units stays bounded below 1 in absolute value as the function is evaluated at a sequence of points which increase to infinity. The Cram\'{e}r condition fails when the random variable has a support that consists of a finite number of points, and hence does not apply to resampling-based distributions such as bootstrap.\footnote{Despite the failure of the Cram\'{e}r condition, the Edgeworth expansion for lattice distributions is well known. See \cite{Bhattacharya/RangaRao:10:NormAppprox}, Chapter 5. See \cite{Kolassa&McCullagh:1990:AS} for a general result of Edgeworth expansion for lattice distributions. \cite{Booth/Hall/Wood:94:JMA} provide the Edgeworth expansion for discrete yet non-lattice distributions.} A standard approach to deal with this issue is to derive an Edgeworth expansion separately for the bootstrap distribution by using an expansion of the empirical characteristic functions. (e.g. \cite{Singh:81:AS} and \cite{Hall:92:Bootstrap}.)

The main contribution of this paper is to provide a finite sample bound for the remainder term in the Edgeworth expansion for a sum of independent random vectors. Using the finite sample bound, one can immediately obtain a uniform-in-$P$ Edgeworth expansion, where the error bound for the remainder term in the expansion is uniform over a collection of probabilities. A notable feature of the Edgeworth expansion is that it admits random vectors which can be discrete, non-lattice distributions. From this result, as shown in the paper, a uniform-in-$P$ Edgeworth expansion for various resampling-based discrete distributions follows as a corollary.

To obtain such an expansion, this paper uses a uniform-in-$P$ version of weak Cram\'{e}r conditions introduced by \cite{Angst/Poly:17:EJP} and obtains a finite sample bound for the error term in the Edgeworth expansion by modifying the proofs of Theorems 20.1 and 20.6 of \cite{Bhattacharya/RangaRao:10:NormAppprox} and Theorem 4.3 of \cite{Angst/Poly:17:EJP}. This paper's finite sample bound reveals that we obtain a uniform-in-$P$ Edgeworth expansion, whenever we have a uniform-in-$P$ bound for the moment of the same order in the Edgeworth expansion.

A uniform-in-$P$ asymptotic approximation is naturally required in a testing set-up with a composite null hypothesis. By definition, a  composite null hypothesis involves a collection of probabilities in the null hypothesis, and the size of a test in this case is its maximal rejection probability over all the probabilities admitted in the null hypothesis. Asymptotic control of size requires uniform-in-$P$ asymptotic approximation of the test statistic's distribution under the null hypothesis. One can apply the same notion to the coverage probability control of confidence intervals as well.

As for uniform-in-$P$ Gaussian approximation, one can obtain the result immediately from a Berry-Esseen bound once appropriate moments are bounded uniformly in $P$. It is worth noting that uniform-in-$P$ Gaussian approximation of empirical processes was studied by \cite{Gine/Zinn:91:AS} and \cite{Sheehy/Wellner:92:AS}. There has been a growing interest in uniform-in-$P$ inference in various nonstandard set-ups in the literature of econometrics in connection with the finite sample stability of inference. (See \cite{Mikusheva:2007:Eca}, \cite{Linton/Song/Whang:2010:JOE}, and \cite{Andrews/Shi:13:Eca}, among others.)

When a test based on a resampling procedure exhibits higher order asymptotic refinement properties, the uniform-in-$P$ Edgeworth expansion can be used to establish a higher order asymptotic size control for the test.   A related work is found in \cite{Hall/Jing:95:AS} who used the uniform-in-$P$ Edgeworth expansion to study the asymptotic behavior of the confidence intervals based on a studentized t statistic. They used a certain smoothness condition for the distributions of the random vectors which excludes resampling-based distributions.

\section{A Uniform-in-$P$ Edgeworth Expansion}
\subsection{Uniform-in-$P$ Weak Cram\'{e}r Conditions}

\cite{Angst/Poly:17:EJP} (hereafter, AP) introduced what they called a weak Cram\'{e}r condition and a mean weak Cram\'{e}r condition which are weaker than the classical Cram\'{e}r condition. Let us prepare notation. Let $\|\cdot\|$ be the Euclidean norm in $\mathbf{R}^d$, i.e., $\|a\|^2 = \text{tr}(a'a)$. The following definition modifies the weak Cram\'{e}r condition introduced by AP into a condition for a collection of probabilities.

\begin{definition}
	\label{def: weak Cramer cond}
	(i) Given $b,c,R>0$, a collection of the distributions $\mathscr{P}$ of a random vector $W$ taking values in $\mathbf{R}^d$ and having characteristic function $\phi_{W,P}$ under $P \in \mathscr{P}$ is said to satisfy the \textit{weak Cram\'{e}r condition} with parameter $(b,c,R)$, if for all $t \in \mathbf{R}^d$ with $\|t\| > R$,
	\begin{eqnarray}
	\label{bound weak Cramer}
		\sup_{P \in \mathscr{P}}|\phi_{W,P}(t)| \le 1 - \frac{c}{\|t\|^b}.
	\end{eqnarray}
	
	(ii) Given $b,c,R>0$, a collection of the joint distributions $\mathscr{P}_n$ of a triangular array of random vectors  $(W_{i,n})_{i=1}^n$ with each $W_{i,n}$ taking values in $\mathbf{R}^d$ and having characteristic function $\phi_{W_{i,n},P}$ under $P \in \mathscr{P}_n$ is said to satisfy the \textit{mean weak Cram\'{e}r condition} with parameter $(b,c,R)$, if for all $t \in \mathbf{R}^d$ with $\|t\| > R$,
	\begin{eqnarray}
	\label{bound mean weak Cramer}		
		\sup_{P \in \mathscr{P}_n}\frac{1}{n}\sum_{i=1}^n|\phi_{W_{i,n},P}(t)| \le 1 - \frac{c}{\|t\|^b}.
	\end{eqnarray}
\end{definition}
The definition of the weak Cram\'{e}r condition modifies that of AP, making the constants $b,c,R$ independent of the probability $P$. This separation between $(b,c,R)$ and $P$ is crucial for our purpose because the probability $P$ in the context of resampling is typically the empirical measure of the data which depends on the sample size. Since the remainder term in the Edgeworth expansion under the weak Cram\'{e}r condition depends on the constants $(b,c,R)$, it is convenient to make these constants independent of the sample size.

\subsection{A Uniform-in-$P$ Edgeworth Expansion under Weak Cram\'{e}r Conditions}

In this section, we present the main result that gives a finite sample bound for the error term in the Edgeworth expansion. Let us prepare notation first. Given each multi-index $\nu = (\nu_1,...,\nu_d)$, with $\nu_k$ being a nonnegative integer, we let $\bar \chi_{\nu,P}$ be the average of the $\nu$-th cumulant of $V_{n,P}^{-1/2}X_{i,n}$ (over $i=1,...,n$). Define for $z = (z_1,...,z_d) \in \mathbf{R}^d$ and a positive integer $j \ge 1$, 
\begin{eqnarray*}
	\bar \chi_{j,P} (z) = j! \sum_{|\nu| = j}\frac{\bar \chi_{\nu,P}}{\nu_1! \dots \nu_d!} z_1^{\nu_1} \dots z_d^{\nu_d}.  
\end{eqnarray*}
For each $j=1,2,...,$ let
\begin{eqnarray}
\label{tilde Pj}
	\tilde P_j(z:\{\bar \chi_{\nu,P}\}) \equiv \sum_{m=1}^s \frac{1}{m!}\left\{\sum_{j_1,...,j_m}^* \frac{\bar \chi_{j_1+2,P}(z)}{(j_1+2)!} \frac{\bar \chi_{j_2+2,P}(z)}{(j_2+2)!} \dots \frac{\bar \chi_{j_m+2,P}(z)}{(j_m+2)!} \right\},
\end{eqnarray}
where $\sum_{j_1,...,j_m}^*$ is the summation over all $m$-tuples $(j_1,...,j_m)$ of positive integers such that $\sum_{k=1}^m j_k = s$. (See (7.3) of \cite{Bhattacharya/RangaRao:10:NormAppprox} (BR, hereafter).) This polynomial has degree $3j$, the smallest order of the terms in the polynomial is $j+2$, and the coefficients in the polynomial involve only $\bar \chi_{\nu,P}$'s with $|\nu| \le j+2$. (Lemma 7.1 of BR, p.52.) Following the convention, we define the derivative operators as follows:
\begin{eqnarray*}
	D_k \equiv \frac{\partial}{\partial x_k}, k=1,...,d, \text{ and } D^\alpha \equiv D_1^{\alpha_1} \cdot \cdot \cdot D_d^{\alpha_d},
\end{eqnarray*} 
for $\alpha = (\alpha_1,...,\alpha_d) \in \{0,1,2,...\}^d$ and $D=(D_1,...,D_d)$. For each $P \in \mathscr{P}_n$, let $Q_{n,P}$ be the distribution of
\begin{eqnarray*}
	\zeta_{n,P} \equiv \frac{1}{\sqrt{n}}\sum_{i=1}^n V_{n,P}^{-1/2}X_{i,n}
\end{eqnarray*}
and define a signed measure $\tilde Q_{n,s,P}$ as follows: for any Borel set $A \subset \mathbf{R}^d$,
\begin{eqnarray*}
	\tilde Q_{n,s,P}(A) \equiv \int_A \sum_{j=0}^{s-2} n^{-{j/2}} \tilde P_j(-D:\{\bar \chi_\nu\})\phi(x)dx,
\end{eqnarray*}
where $\phi$ is the density of the standard normal distribution on $\mathbf{R}^d$. For each $s \ge 1$, define
\begin{eqnarray*}
\rho_{n,s,P} \equiv \frac{1}{n}\sum_{i=1}^n \mathbf{E}_P\|X_{i,n}\|^s.
\end{eqnarray*}
We introduce notation for modulus of continuity for functions: for any Borel measurable real function $f$ on $\mathbf{R}^d$, and any measure $\mu$, we define for $\varepsilon>0$ and $x \in \mathbf{R}^d$,
\begin{eqnarray*}
	\omega_f(x;\varepsilon) \equiv \sup \left\{|f(z)-f(y)|: z,y \in B(x;\varepsilon)\right\}, \text{ and } \bar \omega_f(\varepsilon;\mu) \equiv \int \omega_f(x;\varepsilon)\mu(dx), 
\end{eqnarray*}
where $B(x;\varepsilon)$ is the $\varepsilon$-open ball in $\mathbf{R}^d$ centered at $x$. For any measurable function $f$ and for $s>0$, we define
\begin{eqnarray*}
	M_s(f) \equiv \sup_{x \in \mathbf{R}^d} \frac{|f(x)|}{1 + ||x||^s}.
\end{eqnarray*}
The theorem below is the main result of this paper which is a modification of Theorem 4.3 of AP with the bound made explicit in finite samples. The proof follows the arguments in the proof of Theorem 20.1 of BR and Theorem 4.3 of AP, making bounds in the error explicit in finite samples in each step.

\begin{theorem}
	\label{thm: Edgeworth}
	Suppose that for each $P \in \mathscr{P}_n$, $V_{n,P}$ is positive definite, and the following conditions hold for some integer $s \ge 3$, and $\bar \rho, b,c,R>0$ such that $0 < b < 2/\max\{s-3,1\}$, $c>0$ and $R \ge 1/(16 \bar \rho)$.\medskip
	
    (i) $\sup_{n \ge 1} \sup_{P \in \mathscr{P}_n} \rho_{s,n,P} \le \bar \rho < \infty$.

	(ii) $\mathscr{P}_n$ satisfies the mean weak Cram\'{e}r condition with parameter $(b,c,R)$.\medskip

	Then, there exist constants $n_0 \ge 1$ and $C > 0$ such that for any Borel measurable function $f$ on $\mathbf{R}^d$, for all $n \ge n_0$, and $P \in \mathscr{P}_n$, we have
	\begin{eqnarray*}
		\left| \int f d(Q_{n,P} - \tilde Q_{n,s,P}) \right| \le C n^{-(s-2)/2}\left( 1 + M_s(f)\right) + C \bar \omega_f(n^{-(s-2)/2};\Phi),
	\end{eqnarray*} 
where $\Phi$ denotes the distribution of $N(0,I_d)$, and $n_0$ and $C$ depend only on $s$, $d$, $\bar \rho$, and $(b,c,R)$. 
\end{theorem} 
\medskip

The theorem reveals that under the mean weak Cram\'{e}r condition, the uniform-in-$P$ Edgeworth expansion of $Q_{n,P}$ is essentially obtained by the moment condition that is uniform in $P$ as in Condition (i) of Theorem \ref{thm: Edgeworth}. 

It is worth noting that the bound in Theorem \ref{thm: Edgeworth} depends on $f$ only through $M_s(f)$ and $\bar \omega_f(n^{-(s-2)/2};\Phi)$. When $f$ is such that $|f| \le 1$, we have $M_s(f) \le 1$ and the bound depends on $f$ only through $\bar \omega_f(n^{-(s-2)/2};\Phi)$. In particular, when we take $f$ to be the indicator function of convex subsets of $\mathbf{R}^d$, we obtain the following corollary which is a version of Corollary 20.5 of BR (p.215) with the finite sample bound made explicit here. 

\begin{corollary}
	\label{cor: convex sets}
	Suppose that the conditions of Theorem \ref{thm: Edgeworth} hold. Then, there exist constants $n_0 \ge 1$ and $C > 0$ such that for all $n \ge n_0$, $P \in \mathscr{P}_n$, and for all convex $A \subset \mathbf{R}^d$,
	\begin{eqnarray*}
		\left| Q_{n,P}(A) - \tilde Q_{n,s,P}(A)\right| \le C n^{-(s-2)/2},
	\end{eqnarray*}
    where $n_0$ and $C$ depend only on $s$, $d$, $\bar \rho$, and $(b,c,R)$.
\end{corollary}

The result follows because for any indicator $f(x) = 1\{x \in A\}$ on a convex set $A$, we have $\bar \omega_f(\varepsilon;\Phi) \le c(s,d) \varepsilon$ for any $\varepsilon>0$, where $c(s,d)>0$ is a constant that depends only on $s$ and $d$. (See Corollary 3.2 of BR, p.24.)

\section{An Edgeworth Expansion of Resampling Distributions}

\subsection{Preliminary Results on Weak Cram\'{e}r Condition}
As noted by AP, the weak Cram\'{e}r condition is useful for dealing with distributions obtained from a resampling procedure. To clarify this in our context, let us introduce some notation. For any integers $d,p \ge 1$, a number $h>0$, and any given $d \times p$ matrix $u$ whose $j$-th column is given by $d$-dimensional vector $u_j$, $j=1,...,p$, let
\begin{eqnarray*}
	\mathscr{M}_{d,p}(u;h) \equiv \left\{\sum_{j=1}^p c_j \delta_{u_j}: (c_j)_{j=1}^p \in (h,\infty)^p, \sum_{j=1}^p c_j = 1\right\},
\end{eqnarray*}
where $\delta_{u_j}$ denotes the Dirac measure at $u_j$. If we restrict $c_j = 1/n$ and $p = n$, the class $\mathscr{M}_{d,p}(u;h)$ becomes a singleton that contains
\begin{eqnarray}
\label{resampling distr}
\frac{1}{n}\sum_{j=1}^n  \delta_{u_j}.
\end{eqnarray}
This is the resampling distribution (equivalently, the (nonparametric) bootstrap distribution or the empirical measure) of $(u_1,...,u_n)$, when $(u_1,...,u_n)$ is a realized value of the data $X_n = (X_{n,1}...,X_{n,n})$.

Let $\mathscr{U}_{d,p}(b,c,R,h)$ be the collection of $u$'s such that $\mathscr{M}_{d,p}(u;h)$ does not satisfy the weak Cram\'{e}r condition with parameter $(b,c,R)$. The following proposition is due to AP. (See Proposition 2.4 there.)

\begin{proposition}[\cite{Angst/Poly:17:EJP}]
	\label{prop: resampling AP}
	Suppose that $p \ge 3$ and $b > 1/(p-2)$. Then
	\begin{eqnarray*}
		\lambda\left(\bigcap_{c,R,h>0}\mathscr{U}_{d,p}(2b,c,R,h)\right) = 0,
	\end{eqnarray*}
	where $\lambda$ is Lebesgue measure on $\mathbf{R}^p$.
\end{proposition}

Therefore, the weak Cram\'{e}r condition with parameter $(2b,c,R)$ is generically satisfied by $\mathscr{M}_{d,p}(u;h)$ for some $(c,R,h)$ for Lebesgue almost all $u$'s. However, it is not a trivial task to apply this result to obtain a higher order refinement property of a resampling procedure (even if the data $X_1,...,X_n$ are absolutely continuous with respect to Lebesgue measure.) The main reason is that the error bound in the Edgeworth expansion of such a distribution involves constants, $b,c,R$, which appear in the definition of the weak Cram\'{e}r condition. The constants, $c,R,h$, which make the Lebesgure measure in Proposition \ref{prop: resampling AP}  small may depend on the dimension $p$. In many resampling procedures, this $p$ grows to infinity as the size $n$ of the original sample does so. Therefore, in order to determine the asymptotic behavior of the error bound in the Edgeworth expansion, we need to make explicit how the Lebegue measure in Proposition \ref{prop: resampling AP} behaves as we change $c,R,h$ and $p$.

Instead of working with Lebesgue measure, we directly work with the distribution of the data $X_n=(X_{1,n},...,X_{n,n})$ and, using part of the arguments in the proof of Lemma 5.1 of AP and an exponential tail bound of a U-statistic (Proposition 2.3(b) of \cite{Arcones/Gine:93:AP}), we obtain a finite sample bound for the probability that the empirical measure of $X_n$ does not satisfy the weak Cram\'{e}r condition with parameters $(b,c,R)$. Later we use this result to obtain an Edgeworth expansion for a resampling distribution. This approach, as opposed to one that uses a result like Proposition \ref{prop: resampling AP}, has the additional merit of not having to require $X_n$ to be continuous.

Let $\mathscr{U}_n^*(b,c,R)$ be the collection of $(u_j)_{j=1}^n$, $u_j \in \mathbf{R}^d$, such that the resampling distribution (\ref{resampling distr}) does not satisfy the weak Cram\'{e}r condition with parameter $(b,c,R)$. Proposition \ref{prop: leb bd} below gives an explicit finite sample bound for the probability that an i.i.d. random vector $X_n = (X_{i,n})_{i=1}^n$ does not satisfy the weak Cram\'{e}r condition.

\begin{proposition}
	\label{prop: leb bd}
	Let $X_n = (X_{i,n})_{i=1}^n$ be i.i.d. across $i$'s under $P$ for any $P \in \mathscr{P}_n$. Suppose that there exists a constant $c_R>0$ such that for all $n \ge 1$ and for all $P \in \mathscr{P}_n$,
	\begin{eqnarray}
	\label{bound22}
		c_R \le \sup_{t \in \mathbf{R}^d:\|t\|>R} \frac{1}{2\pi^2} \mathbf{E}_P\left[\inf_{q \in \mathbb{Z}}\left( t' (X_{1,n} - X_{2,n}) - 2 \pi q\right)^2\right],
	\end{eqnarray}
	where $\mathbb{Z} \equiv \left\{...,-1,0,1,...\right\}$.
	
	Then, for all $n \ge 1$,
	\begin{eqnarray*}
		\sup_{P \in \mathscr{P}_n} P\left\{X_n \in \mathscr{U}_n^*\left(b,c_R R^b,R\right)\right\} \le \exp\left( - \frac{c_R^2 n }{2}\right).
	\end{eqnarray*}
\end{proposition}
\medskip

The existence of constant $c_R$ in Proposition \ref{prop: leb bd} can be verified with lower level conditions for $X_n = (X_{i,n})_{i=1}^n$. Since we can always choose $t \in \mathbf{R}^d$ in (\ref{bound22}) such that all the entries except for the first one are zero, let us focus on the case with $d=1$ and note that
\begin{eqnarray*}
	\mathbf{E}_P\left[\inf_{q \in \mathbb{Z}}\left( t(X_{1,n} - X_{2,n}) - 2 \pi q\right)^2\right]
	= t^2 \sum_{k \in \mathbb{Z}} c_P(k;r) \ge R^2 \sum_{k \in \mathbb{Z}} c_P(k;r),
\end{eqnarray*}
where $r \equiv \pi/t$, and
\begin{eqnarray*}
	c_P(k;r) \equiv \left\{\begin{array}{ll}
		\mathbf{E}_P\left[\left((X_{1,n} - X_{2,n}) - r k\right)^21\{X_{1,n} - X_{2,n} \in (rk, r (k+1))]\}\right], & \text{ if } k \text{ is even}\\
		 \mathbf{E}_P\left[\left((X_{1,n} - X_{2,n}) - r(k+1)\right)^21\{X_{1,n} - X_{2,n} \in (rk, r (k+1)]\}\right], & \text{ if } k \text{ is odd}.
		\end{array}
	\right.
\end{eqnarray*}
Suppose that there exist $r \in (0,\pi/R)$ and an integer $k \in \mathbb{Z}$ such that
\begin{eqnarray}
\label{bd234}
	\inf_{n \ge 1} \inf_{P \in \mathscr{P}_n} c_P(k;r) >0.
\end{eqnarray}
Then, the condition (\ref{bound22}) of Proposition \ref{prop: leb bd} holds. This latter condition is very weak. Below, we provide some sufficient conditions.\medskip

\begin{proposition}
	\label{prop: suff}
	Suppose that $\mathscr{P}_n$ in Proposition \ref{prop: leb bd} satisfies that there exist numbers $c_L < c_U$ and $\varepsilon>0$ and $j\in \{1,...,d\}$ such that the $j$-th entry $X_{ij,n}$ of $X_{i,n}$ satisfies either of the following two conditions.
	
	(i) $P\{X_{ij,n} = c_L\} > \varepsilon$ and  $P\{X_{ij,n} = c_U\} > \varepsilon$, for all $P \in \mathscr{P}_n$.
	
	(ii) $\inf_{x \in [c_L,c_U] } f_{j,P}(x) > \varepsilon$, for all $P \in \mathscr{P}_n$, where $f_{j,P}:\mathbf{R} \rightarrow [0,\infty)$ is a measurable map such that for any Borel $A \subset \mathbf{R}$,
	\begin{eqnarray*}
		P\{X_{ij,n} \in [c_L,c_U] \cap A\} = \int_{[c_L,c_U] \cap A} f_{j,P}(x)dx.
	\end{eqnarray*}
	
	Then there exists $c_R>0$ such that Condition (\ref{bound22}) of Proposition \ref{prop: leb bd} holds for all $n \ge 1$ and all $P \in \mathscr{P}_n$.
\end{proposition}

\subsection{An Edgeworth Expansion of Resampling Distributions}
\label{subsec: bootstrap}
\subsubsection{A General Result}
Let us illustrate how the previous results can be applied to obtain a uniform-in-$P$ Edgeworth expansion of a bootstrap distribution of a sum of independent random variables. The Edgeworth expansion of a bootstrap distribution for the i.i.d. random variables is well known in the literature. (See, for example, \cite{Hall:92:Bootstrap}.) The result in this paper is distinct for two reasons. First, the Edgeworth expansion is uniform in $P$, where $P$ runs over a collection of the distributions of the random variable. Second, the Edgeworth expansion follows directly from the Edgeworth expansion for a sum of i.i.d. random variables, due to the use of the weak Cram\'{e}r condition.

Suppose that we have a triangular array of i.i.d. random vectors $X_{i,n}$ taking values in $\mathbf{R}^d$, following distribution $P \in \mathscr{P}_n$, and the bootstrap sample $\{X_{i,n,b}^*\}_{i=1}^n,b=1,...,B$, (i.e., the i.i.d. draws from the empirical measure of $\{X_{i,n}\}_{i=1}^n$), and let the $\sigma$-field generated by $\{X_{i,n}\}_{i=1}^n$ be $\mathscr{F}_n$. Let
\begin{eqnarray*}
	\hat V_n \equiv \frac{1}{n}\sum_{i=1}^n(X_{i,n} - \overline X_n)(X_{i,n} - \overline X_n)'
\end{eqnarray*}
and the conditional distribution of $\sqrt{n}\hat V_n^{-1/2}(\overline X_{n,b}^* - \overline X_n)$ (defined in (\ref{bootstrap test})) given $\mathscr{F}_n$ be denoted by $Q_{n,\mathscr{F}_n}$, where $\overline X_{n,b}^* \equiv \frac{1}{n}\sum_{i=1}^n X_{i,n,b}^*$. Let $\chi_\nu^*$ be the $\nu$-th cumulant of the conditional distribution of $X_{i,n,b}^*$ given $\mathscr{F}_n$. Define for each Borel $A \subset \mathbf{R}^d$,
\begin{eqnarray}
\label{tilde Q}
\tilde Q_{n,\mathscr{F}_n}(A) \equiv \int_A \sum_{j=0}^{s-2} n^{-{j/2}} \tilde P_j(-D:\{\chi_\nu^*\})\phi(x)dx.
\end{eqnarray}
Our purpose is to obtain a finite sample bound for the error term in the approximation of the bootstrap measure $Q_{n,\mathscr{F}_n}$ by its Edgeworth expansion $\tilde Q_{n,\mathscr{F}_n}$. 
Define the event
\begin{eqnarray*}
	\mathscr{E}_{n,0}(\bar \rho) \equiv \left\{\frac{1}{n}\sum_{i=1}^n \|X_{i,n}\|^s \le \bar \rho\right\},
\end{eqnarray*}
and let $\mathscr{E}_{n,1}(c_1)$, $c_1>0$, be the event that the smallest eigenvalue of $\hat V_n$ is bounded from below by $c_1$. Define
\begin{eqnarray*}
	\mathscr{E}_n(\bar \rho,c_1) \equiv \mathscr{E}_{n,0}(\bar \rho) \cap \mathscr{E}_{n,1}(c_1).
\end{eqnarray*}
We are ready to provide an Edgeworth expansion of the bootstrap distribution of $Q_{n,\mathscr{F}_n}$ in the theorem below. The result can be viewed as a variant of the form in (5.29) of \cite{Hall:92:Bootstrap}, p.255, such that the constants are made explicit.\footnote{Note that the constants in (5.29) of \cite{Hall:92:Bootstrap} potentialy depend on $P \in \mathscr{P}_n$, for example, through $C_8$ that appears in (5.16) on page 248.}
\begin{theorem}
	\label{thm: bootstrap}
	Let $s \ge 3$ be a given integer and $\bar \rho>0$ is a number. Suppose that $X_{i,n}$'s are such that (\ref{bound22}) is satisfied for some $c_R>0$ and $0<b< 2/\max\{s-3,1\}$, and $R \ge 1/(16 \bar \rho)$, and let $\mathscr{A}$ be a given subcollection of the Borel $\sigma$-field of $\mathbf{R}^d$. 
	
	Then, there exist constants $n_0 \ge 1$ and $C > 0$ such that for all $n \ge n_0$, and for all $P \in \mathscr{P}_n$,
	\begin{eqnarray}
	\label{sup2} \quad \quad
	&& P\left(\left\{\sup_{A \in \mathscr{A}}\left|Q_{n,\mathscr{F}_n}(A) - \tilde Q_{n,\mathscr{F}_n}(A)\right| - C\bar \omega_{1\{ \cdot \in A\}}\left(n^{-\frac{s-2}{2}};\Phi \right)\le 2 C n^{-\frac{s-2}{2}} \right\} \cap \mathscr{E}_n(\bar \rho,c_1)\right)\\ \notag
	&& \ge P(\mathscr{E}_n(\bar \rho,c_1)) - \exp\left( - \frac{c_R^2 n}{2}\right),
	\end{eqnarray}
	where $C$ and $n_0$ depend only on  $s$, $d$, $\bar \rho$, $b$, $c_R$, and $R$.
	
	Furthermore, if each set $A \in \mathscr{A}$ is convex, then (\ref{sup2}) holds with $\bar \omega_{1\{ \cdot \in A\}}(n^{-(s-2)/2};\Phi)$ replaced by $c(s,d)n^{-\frac{s-2}{2}}$, where $c(s,d)>0$ is a constant that depends only on $s$ and $d$. 
\end{theorem}

In applying the theorem, it remains for us to determine $P(\mathscr{E}_n^c(\bar \rho,c_1))$ and the last supremum in (\ref{sup2}). Following the arguments on pages 245 - 246 of \cite{Hall:92:Bootstrap}, it is not hard to show that $\sup_{P \in \mathscr{P}_n} P(\mathscr{E}_n^c(\bar \rho,c_1)) \le C n^{-\lambda}$ for some constant $C>0$ and $\lambda > (s-2)/2$ under an appropriate uniform-in-$P$ moment condition for $X_{i,n}$.

\subsubsection{A Modulus of Continuity of Resampling Distributions}
Our next result is related to the modulus of continuity of the resampling distribution $Q_{n,\mathscr{F}_n}$. As the distribution $Q_{n,\mathscr{F}_n}$ is discrete, the modulus of continuity around convex sets can be established only up to $n^{-1/2}$ when we use a normal approximation of $Q_{n,\mathscr{F}_n}$. However, using a higher order approximation of $Q_{n,\mathscr{F}_n}$ by $\tilde Q_{n,\mathscr{F}_n}$, we can achieve the modulus of continuity up to the order $n^{-(s-2)/2}$, depending on the moment conditions.

Given $c_2 >0$, define the event
\begin{eqnarray*}
	\mathscr{E}_{n,2}(c_2) &\equiv& \left\{\max_{v\equiv (v_1,...,v_d) \in \mathbb{Z}_+: |v| \le s} \frac{1}{n}\sum_{i=1}^n \prod_{k=1}^d |X_{ik,n}|^{v_k} \le c_2\right\},
\end{eqnarray*}
where $\mathbb{Z}_+ \equiv \{0,1,2,...\}$ and $X_{ik,n}$ denotes the $k$-th entry of $X_{i,n}$. Let
\begin{eqnarray*}
	\mathscr{E}_n(\bar \rho,c_1,c_2) \equiv \mathscr{E}_n(\bar \rho,c_1) \cap \mathscr{E}_{n,2}(c_2).
\end{eqnarray*}
Under appropriate moment conditions, one can show that $P \mathscr{E}_{n,2}^c(c_2) = O(n^{-\lambda})$ as $n \rightarrow \infty$, for some $\lambda>\delta$ and large $c_2 >0$. (See Chapter 5 of \cite{Hall:92:Bootstrap} for details.)

\begin{corollary}
	\label{cor: mod cont}
	Suppose that the assumptions of Theorem \ref{thm: bootstrap} hold, and let $C$ and $n_0$ be as in Theorem \ref{thm: bootstrap}. For any set $A \subset \mathbf{R}^d$ and for $\eta>0$, $A^\eta$ be its $\eta$-enlargement, i.e., $A^\eta = \{x \in \mathbf{R}^d: \exists y \in A, \text{ s.t. } \|x - y\| \le \eta\}$. Let $\mathscr{C}_d$ be the collection of convex sets.
	
	Then for all $n \ge n_0$, $\eta>0$, and $P \in \mathscr{P}_n$,
	\begin{eqnarray*}
		&& P\left(\left\{\sup_{A \in \mathscr{C}_d}\left| Q_{n,\mathscr{F}_n}(A^\eta) - Q_{n,\mathscr{F}_n}(A) \right| \le C'n^{-(s-2)/2} + c_1(s,d) c_2 \eta\right\}\cap \mathscr{E}_n(\bar \rho,c_1,c_2)\right)\\
		&& \ge P(\mathscr{E}_n(\bar \rho,c_1,c_2)) - \exp\left(-\frac{c_R^2 n}{2} \right),
	\end{eqnarray*}
	where $C'>0$ is a constant that depends only on $C$, $s$ and $d$, and $c_1(s,d)$ is a constant that depends only on $s$ and $d$.
\end{corollary}

\subsubsection{Uniform-in-$P$ Edgeworth Expansion of a t-Test Statistic} Suppose that $\{W_{i,n}\}_{i=1}^n$ is a triangular array of random variables which are i.i.d. drawn from a common distribution $P$. Let us assume that this distribution belongs to a collection $\mathscr{P}_n$ of distributions. Let $\{W_{i,n,b}^*\}_{i=1}^n, b=1,...,B$, be the bootstrap sample drawn with replacement from the empirical distribution of $\{W_{i,n}\}_{i=1}^n$. Define the sample variance:
\begin{eqnarray*}
	s_{n,b}^{*2} \equiv \frac{1}{n}\sum_{i=1}^n(W_{i,n,b}^* - \overline W_{n,b}^*)^2,
\end{eqnarray*}
where $\overline W_{n,b}^* \equiv \frac{1}{n}\sum_{i=1}^n W_{i,n,b}^*$. Then, we are interested in the uniform-in-$P$ Edgeworth expansion of the bootstrap distribution of the following:
\begin{eqnarray*}
	\frac{1}{s_{n,b}^*\sqrt{n}}\sum_{i=1}^n (W_{i,n,b}^*-\overline W_n),
\end{eqnarray*}
where $\overline W_n \equiv \frac{1}{n}\sum_{i=1}^n W_{i,n}$. Following Chapter 5 of \cite{Hall:92:Bootstrap}, we let
\begin{eqnarray}
\label{spec}
	\overline X_{n,b}^* \equiv \left(\overline W_{n,b}^*, \frac{1}{n}\sum_{i=1}^n W_{i,n,b}^{*2}\right),
\end{eqnarray}
and write
\begin{eqnarray}
\label{bootstrap test2}
\frac{1}{s_{n,b}^*\sqrt{n}}\sum_{i=1}^n (W_{i,n,b}^*-\overline W_n) =  \sqrt{n} g_n(\overline X_{n,b}^*),
\end{eqnarray}
where $g_n(x) =  (x_1 - \overline W_n)/\sqrt{x_2 - x_1^2}$, $x = (x_1,x_2)$, with $x_2 > x_1^2$. Our focus is on the Edgeworth expansion of the bootstrap distribution of the test statistic of the form:
\begin{eqnarray}
\label{bootstrap test}
		T_{n,b}^* \equiv \sqrt{n} g_n(\overline X_{n,b}^*).
\end{eqnarray}

Let $X_{i,n} = (W_{i,n},W_{i,n}^2)$ and $\mathscr{F}_n$ be the $\sigma$-field of $X_n = (X_{i,n})_{i=1}^n$. Let $Q_{n,\mathscr{F}_n}$ be the conditional distribution of $\sqrt{n} \hat V_n^{-1/2}(\overline X_{n,b}^* - \overline X_n)$ given $\mathscr{F}_n$. Define for any $t \in \mathbf{R}$,
\begin{eqnarray*}
	\hat f_{n,t}(x) \equiv 1\left\{\sqrt{n}g_n\left(\overline X_n + \frac{\hat V_n^{1/2} x}{\sqrt{n}}\right) \le t \right\}.
\end{eqnarray*}
Then the map $Q_{n,\mathscr{F}_n}(\hat f_{n,t})$ represents the bootstrap distribution of $T_{n,b}^*$. Note that we cannot use the convex set approach to deal with the last term in (\ref{sup2}) as before, because the set of $x$'s such that $\hat f_{n,t}(x) = 1$ is not convex. Instead, we invoke Lemma 5.3 of \cite{Hall:92:Bootstrap}.

We define for constants $s,c_3>0$,
\begin{eqnarray*}
	\mathscr{E}_{n,3}(c_3) &\equiv& \left\{\max_{|\alpha|\le s + 3}\left|D^\alpha g_n(\overline X_n)\right| \le c_3\right\},
\end{eqnarray*}
and denote by $\mathscr{E}_{n,3}(c_3)$, $c_3>0$, the event that the largest eigenvalue of $\hat V_n$ is bounded by $c_3$. Define
\begin{eqnarray*}
	\mathscr{E}_{n}(\bar \rho,c_1,c_2,c_3) \equiv \mathscr{E}_{n,0}(\bar \rho) \cap \mathscr{E}_{n,1}(c_1) \cap \mathscr{E}_{n,2}(c_2) \cap \mathscr{E}_{n,3}(c_3).
\end{eqnarray*}
Using Lemma 5.3 of \cite{Hall:92:Bootstrap}, we can obtain the following corollary.

\begin{corollary}
	\label{cor: bootstrap}
	Suppose that the conditions of Theorem \ref{thm: bootstrap} hold. Then, there exist constants $n_0' \ge 1$ and $C > 0$ such that for all $n \ge n_0'$ and all $P \in \mathscr{P}_n$,
	\begin{eqnarray}
	\label{sup3} \quad \quad
	&& P\left(\left\{\sup_{t \in \mathbf{R}} |Q_{n,\mathscr{F}_n}(\hat f_{n,t}) - \tilde Q_{n,\mathscr{F}_n}(\hat f_{n,t})| \le C n^{-\frac{s-2}{2}} \right\} \cap \mathscr{E}_{n}(\bar \rho,c_1,c_2,c_3) \right)\\ \notag
	&& \ge P\left(\mathscr{E}_{n}(\bar \rho,c_1,c_2,c_3)\right) - \exp\left(-\frac{c_R^2 n}{2} \right),
	\end{eqnarray}
	where $C$ and $n_0'$ depend only on $s$, $d$, $\bar \rho$, $b$, $c_R$, $R$, $n_0$, and $c_1,c_2,c_3$.
\end{corollary}

As for the probability $P (\mathscr{E}_{n}^c(\bar \rho,c_1,c_2,c_3))$, we can find its finite sample bound using standard arguments, and show it to be $o(n^{-(s-2)/2})$ uniformly over $P \in \mathscr{P}_n$, when there is a uniform bound for the population moments, and a uniform lower and upper bound for the eigenvalues of $V_{n,P} \equiv \mathbf{E}_P[(X_{i,n} - \mathbf{E}_PX_{i,n})(X_{i,n} - \mathbf{E}_PX_{i,n})']$. Details can be furnished following the same arguments in Chapter 5 of \cite{Hall:92:Bootstrap}.

\section{Proofs}
We introduce notation. For a measure or a signed measure $\mu$ and a measurable function $f$, we write
\begin{eqnarray*}
	\mu(f) \equiv \int f d\mu,
\end{eqnarray*}
for brevity. Recall that the Fourier-Stieltjes transform of a signed measure $\mu$ on $\mathbf{R}^d$ is a complex-valued function on $\mathbf{R}^d$:
\begin{eqnarray*}
	\hat \mu(t) \equiv \int_{\mathbf{R}^d} e^{\text{i}t'x} \mu(dx), t \in \mathbf{R}^d,
\end{eqnarray*}
where $\text{i} = \sqrt{-1}$. Recall that $\|\cdot\|$ denotes the Euclidean norm on $\mathbf{R}^d$, i.e., $\|a\|^2 = a'a$, $a \in \mathbf{R}$. We also define $|a| \equiv \sum_{k=1}^d|a_k|$, for any $a = (a_1,...,a_d)' \in \mathbf{R}^d$.

Define
\begin{eqnarray}
\label{YZ}
	Y_{i,n} \equiv X_{i,n}1\{\|X_{i,n}\| \le \sqrt{n}\}, \text{ and } Z_{i,n} \equiv Y_{i,n} - \mathbf{E}_P Y_{i,n},
\end{eqnarray}
and let $\phi_{i,n,P}(t) \equiv \mathbf{E}_P[\exp(\text{i}t'Z_{i,n})]$ and $\phi_{X_{i,n},P}(t) \equiv \mathbf{E}_P[\exp(\text{i}t'X_{i,n})]$. Let us begin with an auxiliary lemma which is a uniform-in-$P$ modification of Proposition 2.10 and Lemma 5.4 of AP.
\begin{lemma}
	\label{lemma: prop 2.10}
	Let $\{X_{i,n}\}_{i=1}^n$ be a triangular array of independent random vectors taking values in $\mathbf{R}^d$, and let $\mathscr{P}_n$ be the collection of the joint distributions of $(X_{i,n})_{i=1}^n$ such that $\mathscr{P}_n$ satisfies the mean weak Cram\'{e}r condition with parameter $(b,c,R)$ for some $b,c,R>0$. Furthermore, assume that
	\begin{eqnarray*}
		\sup_{n \ge 1} \sup_{P \in \mathscr{P}_n} \frac{1}{n}\sum_{i=1}^n \mathbf{E}_P\left[\|X_{i,n}\|\right] < \infty.
	\end{eqnarray*}
	
	Then there exist $\varepsilon>0$ and $n_0 \ge 1$ such that for all $0 < r < R$ and all $n \ge n_0$,
	\begin{eqnarray*}
		\sup_{r \le \| u \| \le R } \sup_{P \in \mathscr{P}_n}\frac{1}{n}\sum_{i=1}^n |\phi_{X_{i,n},P}(u)| \le 1 - \varepsilon,
	\end{eqnarray*}
	where $\varepsilon$ and $n_0$ depend only on $b,c,R$ and $r$.
	
	Suppose further that for some $s \ge 3$ and $\bar \rho>0$,
	\begin{eqnarray*}
		\sup_{n \ge 1} \sup_{P \in \mathscr{P}_n} \frac{1}{n}\sum_{i=1}^n \mathbf{E}_P\left[\|X_{i,n}\|^s\right] \le \bar \rho.
	\end{eqnarray*}

	Then, there exist $\varepsilon>0$ and $n_0$ such that for all $0 < r < R$ and all $n \ge n_0$,
	\begin{eqnarray*}
		\sup_{r  \le \|t\| \le R}\sup_{P \in \mathscr{P}_n}\frac{1}{n}\sum_{i=1}^n |\phi_{i,n,P}(t)| \le 1 - \varepsilon,
	\end{eqnarray*} 
and for all $t \in \mathbf{R}^d$ such that $R < \|t\| <(c/(4 \bar \rho))^{1/b} n^{s/(2b)}$, 
\begin{eqnarray*}
	\sup_{P \in \mathscr{P}_n}\frac{1}{n}\sum_{i=1}^n \left|\phi_{i,n,P}(t) \right| \le 1 - \frac{c}{2\|t\|^b},
\end{eqnarray*}
where $\varepsilon>0$ and $n_0$ depend only on $s$, $\bar \rho$, $b$, $c$, $R$ and $r$.
\end{lemma}
\medskip

\noindent \textbf{Proof:} The first statement is obtained by following the proof of Proposition 2.10 of AP. More specifically, suppose to the contrary that we have
\begin{eqnarray*}
	\limsup_{n \rightarrow \infty} \sup_{r \le \| u \| \le R } \sup_{P \in \mathscr{P}_n}\frac{1}{n}\sum_{i=1}^n |\phi_{X_{i,n},P}(u)| = 1.
\end{eqnarray*}
Then by the definition of supremum, there exists a sequence $(u_n,P_n)$ such that $r \le \|u_n\| \le R$ and $P_n \in \mathscr{P}_n$, and
\begin{eqnarray*}
	\limsup_{n \rightarrow \infty} \sup_{r \le \| u \| \le R } \sup_{P \in \mathscr{P}_n}\frac{1}{n}\sum_{i=1}^n |\phi_{X_{i,n},P}(u)|  = \limsup_{n \rightarrow \infty} \frac{1}{n}\sum_{i=1}^n |\phi_{X_{i,n},P_n}(u_n)|.
\end{eqnarray*}
We can choose a subsequence $\{n(k)\}_{k=1}^\infty \subset \{n\}_{n \ge 1}$ such that
\begin{eqnarray*}
		\limsup_{n \rightarrow \infty} \frac{1}{n}\sum_{i=1}^n |\phi_{X_{i,n},P_n}(u_n)| = \lim_{k \rightarrow \infty} \frac{1}{n(k)}\sum_{i=1}^{n(k)} |\phi_{X_{i,n(k)},P_{n(k)}}(u_{n(k)})|.
\end{eqnarray*}
We can follow the rest of the proof in AP, p.6-7, to arrive at a contradiction to the weak Cram\'{e}r condition.

As for the second result, we take supremum over $P \in \mathscr{P}_n$ on both sides of (5.17) in Lemma 5.4 of AP to obtain that
\begin{eqnarray*}
	\sup_{r  \le \|t\| \le R} \sup_{P \in \mathscr{P}_n}\frac{1}{n}\sum_{i=1}^n \left|\phi_{i,n,P}(t) \right| \le \sup_{r  \le \|t\| \le R} \sup_{P \in \mathscr{P}_n}\frac{1}{n}\sum_{i=1}^n \left|\phi_{i,n,P}(t) \right| + \frac{2 \bar \rho}{n^{s/2}}.
\end{eqnarray*}
The proof can be completed following the same arguments in the proof of the lemma using Definition \ref{def: weak Cramer cond}. $\blacksquare$\medskip

\noindent \textbf{Proof of Theorem \ref{thm: Edgeworth}:} We walk through the steps in the proofs of Theorem 20.1 of BR and Theorem 4.3 of AP, making bounds in the error explicit in finite samples. Throughout the proof, we assume without loss of generality that $V_{n,P} = I_d$. While the proof of Theorem 20.1 and Theorem 4.3 assumes that $X_{i,n}$'s are i.i.d., here we make it explicit that they are allowed to be non-identically distributed. In the proof, as in BR, we use notation $c_j(x)$ to denote a constant that depends only on $x$, so that, for example, $c_1(s,d)$ is a constant that depends only on $s$ and $d$. In following the proof of Theorem 20.1 of BR, we take $s'$ in the theorem to be $s$ for simplicity. 

Recall the definitions of $Y_{i,n}$ and $Z_{i,n}$ in (\ref{YZ}). Let $Q_{n,P}'$ and $Q_{n,P}''$ be the distributions of $\frac{1}{\sqrt{n}}\sum_{i=1}^nZ_{i,n}$ and $\frac{1}{\sqrt{n}}\sum_{i=1}^nY_{i,n}$ respectively. We take
\begin{eqnarray*}
	a_n \equiv \frac{1}{\sqrt{n}}\sum_{i=1}^n \mathbf{E}[Y_{i,n}],
\end{eqnarray*}
and define
\begin{eqnarray}
\label{fan}
	f_{a_n}(x) \equiv f(x + a_n).
\end{eqnarray}
Define for any constant $a>0$ and integer $s \ge 1$,
\begin{eqnarray*}
	\bar \triangle_{n,s,P}(a) \equiv \frac{1}{n}\sum_{i=1}^n \mathbf{E}_P\left[||X_{i,n}||^s 1\left\{||X_{i,n}|| \ge a \sqrt{n}\right\}\right].
\end{eqnarray*}
Let $\delta>0$ be such that
\begin{eqnarray}
\label{bound}
\frac{s-2}{2} < \delta < \frac{1}{b} + \frac{1}{2}.
\end{eqnarray}
Existence of such $\delta>0$ is ensured by the condition that $b< 2/\max\{s-3,1\}$ made in the theorem.
 
Using the identity $Q_{n,P}''(f) = Q_{n,P}'(f_{a_n})$, we bound
\begin{eqnarray*}
	|Q_{n,P}(f) - \tilde Q_{n,s,P}(f)| \le A_{1n} + A_{2n} + A_{3n},
\end{eqnarray*}
where
\begin{eqnarray*}
	A_{1n} &\equiv& |Q_{n,P}(f) - Q_{n,P}''(f)|,\\
	A_{2n} &\equiv& |Q_{n,P}'(f_{a_n}) - \tilde Q_{n,s,P}(f_{a_n})|, \text{ and }\\
	A_{3n} &\equiv& |\tilde Q_{n,s,P}(f_{a_n}) - \tilde Q_{n,s,P}(f)|.
\end{eqnarray*}
By (20.8) of BR, p.208, and (20.12) of BR, p.209, we find that (taking $s' = s$ there)
\begin{eqnarray}
\label{bound1}
A_{1n} + A_{3n} &\le& c_1(s,d) M_s(f) n^{-(s-2)/2} \bar \Delta_{n,s,P}(1) \\ \notag
 &\le& c_1(s,d) M_s(f) n^{-(s-2)/2} \bar \rho.
\end{eqnarray}
Let us focus on finding a bound for $A_{2,n}$. Let $D_{n,P} \equiv \frac{1}{n}\sum_{i=1}^n \text{Var}_P(Z_{i,n})$ and define a signed measure $\tilde Q_{n,s,P}'$ as 
\begin{eqnarray}
\label{tilde Q'}
\tilde Q_{n,s,P}'(A) \equiv \int_A \sum_{j=0}^{s-2} n^{-j/2} \tilde P_j(-D;\{\bar \chi_{\nu,P}\}) \phi_{0,D_{n,P}}(x) dx, \text{ for any Borel set } A,
\end{eqnarray}
where $\phi_{0,D_{n,P}}$ denotes the density of $N(0,D_{n,P})$. We bound
\begin{eqnarray}
\label{bound2}
A_{2n} &\le& |\tilde Q_{n,s,P}'(f_{a_n}) - \tilde Q_{n,s,P}(f_{a_n})| + A_{2n}'\\ \notag
&\le& c_2(s,d) M_s(f) n^{-(s-2)/2} \bar \Delta_{n,s,P}(1) + A_{2n}',	
\end{eqnarray}
where $A_{2n}' \equiv |Q_{n,P}'(f_{a_n}) - \tilde Q_{n,s,P}'(f_{a_n})|$, and the last inequality comes from (20.13) of BR, p.209 (again, taking $s' = s$). Let us focus on $A_{2n}'$. Define
\begin{eqnarray*}
	H_{n,P} \equiv Q_{n,P}' - \tilde Q_{n,s,P}'.
\end{eqnarray*}
For any positive number $0<\varepsilon<1$, we define $K_\varepsilon$ to be a probability measure such that 
\begin{eqnarray*}
	K_\varepsilon(\{x \in \mathbf{R}^d: \|x\| < \varepsilon \}) = 1,
\end{eqnarray*}
and 
\begin{eqnarray*}
	|D^\varepsilon \hat K_\varepsilon(t)| \le c_3(s,d) \varepsilon^{|\alpha|} \exp\left( - (\varepsilon\|t\|)^{1/2} \right),
\end{eqnarray*}
for all $t \in \mathbf{R}^d$ and all $|\alpha| \le s + d + 1.$ Then, by (20.17) of BR, p.210, we find that
\begin{eqnarray}
\label{bound3}
|H_{n,P}(f_{a_n})| &\le& M_s(f) \int \left(1+\left(\|x\| + \varepsilon + \|a\|\right)^s\right)|H_{n,P} * K_\varepsilon|(dx) \\ \notag
&& + \bar \omega_{f_{a_n}}\left( 2 \varepsilon; |\tilde Q_{n,s+d,P}'|\right) \equiv B_{1n} + B_{2n}, \text{ say}.
\end{eqnarray}
From (20.19) of BR, p.210, we bound
\begin{eqnarray*}
	B_{1n} \le c_4(s,d) \max_{0 \le |\beta| \le s + d + 1} \int |D^\beta(\hat H_{n,P}\hat K_\varepsilon)(t)|dt.
\end{eqnarray*}
We write
\begin{eqnarray*}
	D^\beta(\hat H_{n,P} \hat K_\varepsilon) = \sum_{0 \le \alpha \le \beta} c_5(\alpha, \beta)(D^{\beta- \alpha} \hat H_{n,P})(D^\alpha \hat K_\varepsilon).
\end{eqnarray*}
Now, let us focus on finding a bound for
\begin{eqnarray}
\label{int DH DK}
\int\left| (D^{\beta- \alpha} \hat H_{n,P}(t))(D^\alpha \hat K_\varepsilon(t)) \right|dt.
\end{eqnarray}
For any $d\times d$  matrix $A$, we define the operator norm $\|A\| \equiv \sup_{x \in \mathbf{R}^d: \|x\| \le 1} \|A x\|$ as usual. When $A$ is symmetric, we have $\|A^2\| = \|A\|^2$. Define for some constant $c_6(s,d)$ which depends only on $s$ and $d$,
\begin{eqnarray}
\label{An}
A_n \equiv \frac{c_6(s,d) n^{1/2}}{\left(\frac{1}{n}\sum_{i=1}^n \mathbf{E}_P\| D_{n,P}^{-1/2} Z_{i,n}\|^{s+d+1}\right)^{1/(s+d-1)} \Lambda_n^{1/2}},
\end{eqnarray}
where $\Lambda_n$ denotes the maximum eigenvalue of $D_{n,P}$. We bound the integral in (\ref{int DH DK}) by
\begin{eqnarray}
\label{dec}
&& \int_{\|t\| \le A_n}\left| (D^{\beta- \alpha} \hat H_{n,P}(t))(D^\alpha \hat K_\varepsilon(t)) \right|dt\\ \notag
&& + \int_{\|t\| > A_n}\left| (D^{\beta- \alpha} \hat H_{n,P}(t))(D^\alpha \hat K_\varepsilon(t)) \right|dt.
\end{eqnarray}
By (20.21) of BR, p.210, we can choose $c_6(s,d)$ in the definition of $A_n$ so that the leading term in (\ref{dec}) is bounded by
\begin{eqnarray}
\label{bd0}
&& c_7(s,d) n^{-(s+d-1)/2}\frac{1}{n}\sum_{i=1}^n \mathbf{E}_P\left[\|D_{n,P}^{-1/2} Z_{i,n}\|^{s+d+1} \right]\\ \notag
&\le& c_7(s,d) n^{-(s+d-1)/2} \|D_{n,P}^{-1/2}\|^{s+d+1} \frac{1}{n}\sum_{i=1}^n \mathbf{E}_P\left[ \| Z_{i,n}\|^{s+d+1} \right]\\ \notag
&=& c_7(s,d) n^{-(s+d-1)/2} \|D_{n,P}^{-1}\|^{(s+d+1)/2} \frac{1}{n}\sum_{i=1}^n \mathbf{E}_P\left[ \| Z_{i,n}\|^{s+d+1} \right]
\end{eqnarray}
From (14.23) and (14.24) of BR, p.125, and our assumption that $V_{n,P} = I_d$, and because $D_{n,P}$ is symmetric and $\|D_{n,P}\|$ denotes the operator norm of $D_{n,P}$, we find that
\begin{eqnarray}
\label{bd11}
\|D_{n,P} - I_d\| \equiv \sup_{\|t\| \le 1} t'(D_{n,P} - I_d)t &\le& 2d n^{-(s-2)/2} \bar \Delta_{n,s,P}(1)\\ \notag
  &\le& 2d n^{-(s-2)/2} \bar \rho,
\end{eqnarray}
and
\begin{eqnarray}
\label{bd D}
\|D_{n,P}^{-1}\| = \|(I_d + (D_{n,P} - I_d))^{-1}\| &\le& \frac{1}{1 - \|D_{n,P} - I_d\|}\\ \notag
&\le&  \frac{1}{1 - 2d n^{-(s-2)/2} \bar \Delta_{n,s,P}(1)}.
\end{eqnarray}
As noted in (20.23) of BR, p.211,
\begin{eqnarray}
\label{bd1}
\mathbf{E}_P\|Z_{i,n}\|^{s+d+1} \le 2^{s+d+1} \mathbf{E}_P\|Y_{i,n}\|^{s+d+1}.
\end{eqnarray}
For any $\xi \in (0,1]$, we can bound
\begin{eqnarray}
\label{bd2}
\quad \quad \frac{1}{n}\sum_{i=1}^n \mathbf{E}_P\|Y_{i,n}\|^{s+d+1} &\le& (\xi n^{1/2})^{d+1} \frac{1}{n}\sum_{i=1}^n \mathbf{E}_P\left[1\{\|X_{i,n}\| \le \xi n^{1/2}\} \|X_{i,n}\|^s\right]\\ \notag
&& + n^{(d+1)/2} \frac{1}{n}\sum_{i=1}^n \mathbf{E}_P\left[1\{\xi n^{1/2} < \|X_{i,n}\| \le n^{1/2}\} \|X_{i,n}\|^s\right]\\ \notag
&\le& n^{(d+1)/2}\left\{\xi^{d+1} \rho_{n,s,P} + \bar \Delta_{n,s,P}(\xi)\right\} \le 2 n^{(d+1)/2} \bar \rho.
\end{eqnarray}
This gives a bound for $\frac{1}{n}\sum_{i=1}^n\mathbf{E}_P\|Z_{i,n}\|^{s+d+1}$ through (\ref{bd1}). Using this bound and (\ref{bd D}), we find that for any $\xi>0$,
\begin{eqnarray}
\label{ineq}
\frac{1}{n}\sum_{i=1}^n \mathbf{E}_P\left[\|D_{n,P}^{-1/2} Z_{i,n}\|^{s+d+1} \right] \le \frac{2^{s+d+1} n^{(d+1)/2}\left(\xi^{d+1} \rho_{n,s,P} + \bar \Delta_{n,s,P}(\xi) \right)}{\left(1 - 2dn^{-(s-2)/2}\bar \Delta_{n,s,P}(1) \right)^{(s+d+1)/2}},
\end{eqnarray}
and from (\ref{bd0}), we obtain that
\begin{eqnarray}
\label{bound A}
&& \quad \int_{\|t\| \le A_n}\left| (D^{\beta- \alpha} \hat H_{n,P}(t))(D^\alpha \hat K_\varepsilon(t)) \right|dt \\ \notag
&\le&  \frac{c_7'(s,d) n^{-(s-2)/2}\left(\xi^{d+1} \rho_{n,s,P} + \bar \Delta_{n,s,P}(\xi) \right)}{\left(1 - 2dn^{-(s-2)/2}\bar \Delta_{n,s,P}(1) \right)^{(s+d+1)/2}},
\end{eqnarray}
where $c_7'(s,d) = c_7(s,d) 2^{s+d+1}$.

Let us turn to the second integral in (\ref{dec}). Note first that
\begin{eqnarray*}
	\Lambda_n \le \|D_{n,P}\| \le \|D_{n,P} - I_d\| + \|I_d\| &\le& 1 +  2d n^{-(s-2)/2} \bar \Delta_{n,s,P}(1)\\
	&\le& 1 +  2d n^{-(s-2)/2} \bar \rho,
\end{eqnarray*}
by (\ref{bd11}). Using this and noting (\ref{ineq}), we find from the definition of $A_n$ in (\ref{An}) that
\begin{eqnarray}
\label{lower bound An}
A_n &\ge& \frac{c_6(s,d) n^{\frac{s-2}{2(s+d-1)}}}{\left(\xi^{d+1} \rho_{n,s,P} + \bar \Delta_{n,s,P}(\xi)\right)^{\frac{1}{s+d-1}}} \frac{\left(1 - 2dn^{-(s-2)/2}\bar \Delta_{n,s,P}(1) \right)^{\frac{s+d+1}{2(s+d-1)}}}{\left(1 + 2dn^{-(s-2)/2}\bar \Delta_{n,s,P}(1) \right)^{1/2}}\\ \notag
&\ge& \frac{c_6(s,d) n^{\frac{s-2}{2(s+d-1)}}}{\left(2 \bar \rho \right)^{\frac{1}{s+d-1}}} \frac{\left(1 - 2dn^{-(s-2)/2}\bar \rho \right)^{\frac{s+d+1}{2(s+d-1)}}}{\left(1 + 2dn^{-(s-2)/2}\bar \rho \right)^{1/2}}.
\end{eqnarray}
We take 
\begin{eqnarray*}
	q_n \equiv \frac{\sqrt{n}}{16 \bar \rho},
\end{eqnarray*}
and, as in (20.26) of BR, p.211, let us bound
\begin{eqnarray*}
	\int_{\|t\| > A_n}\left| (D^{\beta- \alpha} \hat H_{n,P}(t))(D^\alpha \hat K_\varepsilon(t)) \right|dt \le I_1 + I_2 + I_3,
\end{eqnarray*}
where
\begin{eqnarray*}
	I_1 &\equiv& \int_{\|t\| > q_n}\left| (D^{\beta- \alpha} \hat Q_{n,P}'(t))(D^\alpha \hat K_\varepsilon(t)) \right|dt\\ \notag
	I_2 &\equiv& c_8(s,d) \int_{\|t\| > A_n} \left( 1 + \|t\|^{|\beta - \alpha|}\right) \exp\left( - \frac{5}{24}\|t\|^2\right)dt, \text{ and } \\ \notag
	I_3 &\equiv&  \int_{\|t\| > A_n} \left|D^{\beta - \alpha} \sum_{j=0}^{s+d-2} n^{-j/2}\tilde P_j(\text{i}t;\{\bar \chi_{\nu,P}\}) \exp \left(-\frac{1}{2} t' D_{n,P} t \right) \right|dt.
\end{eqnarray*}
Let us deal with $I_2$ first. Choose any large number $w \ge 1$. We bound
\begin{eqnarray*}
	I_2 \le c_8(s,d) A_n^{-w}\int_{\|t\| > A_n} \left( 1 + \|t\|^{|\beta - \alpha|}\right) \|t\|^w \exp\left( - \frac{5}{24}\|t\|^2\right)dt. 
\end{eqnarray*}
The last integral is bounded for any fixed $w \ge 1$. In the light of the lower bound for $A_n$ in (\ref{lower bound An}), we find that for any $w \ge 1$, there exist a constant $C$ and a number $n_{0,1} \ge 1$ which depend only on $(s,d,\bar \rho,w)$ such that for all $n \ge n_{0,1}$,
\begin{eqnarray*}
	I_2 \le C n^{-\frac{(s-2)w}{2(s+d-1)}}.
\end{eqnarray*}
Since the integrand in $I_3$ involves $\exp(-t'D_{n,P}t/2)$, we obtain the same conclusion for $I_3$.

The rest of the proof focuses on finding a finite sample bound for $I_1$. In doing so, we switch to the proof of Theorem 4.3 of AP. By (5.15) of AP, p.19, we bound
\begin{eqnarray*}
	I_1 \le c_9(s,d) \varepsilon^{|\alpha|} 2^{|\beta - \alpha|} \sum_{\gamma \in (\mathbb{Z}_+^d)^n: |\gamma| = |\beta - \alpha|}\left(\begin{matrix} \beta - \alpha\\ \gamma \end{matrix} \right) J_\gamma(n,\varepsilon),
\end{eqnarray*}
where
\begin{eqnarray*}
	J_\gamma(n,\varepsilon) \equiv n^{d/2}\int_{\|t\| \ge q_n/\sqrt{n}}\left|\prod_{i=1,\gamma_i = 0}^n \phi_{i,n,P}(t) \right| e^{-(\sqrt{n}\varepsilon\|t\|)^{1/2}} dt,
\end{eqnarray*}
and
\begin{eqnarray*}
	\left(\begin{matrix} \beta - \alpha\\ \gamma \end{matrix} \right) = \frac{|\beta - \alpha|!}{\prod_{i=1}^n |\gamma_i|!}.
\end{eqnarray*}
We let
\begin{eqnarray*}
	r\equiv \frac{1}{16 \bar \rho},
\end{eqnarray*}
so that $q_n /\sqrt{n} = r$ for all $n \ge 1$. With this choice of $r>0$, from (5.16) of AP, p.20, we deduce that for all $n \ge n_1$,
\begin{eqnarray*}
	J_\gamma(n,\varepsilon) \le n^{d/2}\exp\left((n - |\gamma|) \log\left(\frac{n}{n - |\gamma|} \right) \right)K_\gamma(n,\varepsilon),
\end{eqnarray*}
where
\begin{eqnarray*}
	K_\gamma(n,\varepsilon) \equiv \int_{\|t\| \ge r} \exp\left((n - |\gamma|) \log\left(\frac{1}{n}\sum_{i=1}^n |\phi_{i,n,P}(t)| \right) \right)e^{-(\varepsilon\sqrt{n}\|t\|)^{1/2}} dt.
\end{eqnarray*}
By Lemma \ref{lemma: prop 2.10}, there exists $\bar n \ge 1$ such that for all $n \ge \bar n$ and all $t \in \mathbf{R}^d$ satisfying $R < \|t\| \le (c/(4\bar \rho))^{1/b} n ^{s/(2b)}$,
\begin{eqnarray*}
	\sup_{P \in \mathscr{P}_n} \frac{1}{n}\sum_{i=1}^n |\phi_{i,n,P}(t)| \le 1 - \frac{c}{2\|t\|^b}.
\end{eqnarray*}
Take $n_0' = \max\{\bar n,n_{0,1}\}$, and write
\begin{eqnarray*}
	K_\gamma(n,\varepsilon) = K_{\gamma}^1(n,\varepsilon) + K_{\gamma}^2(n,\varepsilon) + K_{\gamma}^3(n,\varepsilon),
\end{eqnarray*}
where
\begin{eqnarray*}
	K_\gamma^1(n,\varepsilon) &\equiv& \int_{r \le \|t\| \le R} \exp\left((n - |\gamma|) \log\left(\frac{1}{n}\sum_{i=1}^n |\phi_{i,n,P}(t)| \right) \right)e^{-(\varepsilon\sqrt{n}\|t\|)^{1/2}} dt \\ \notag
	K_\gamma^2(n,\varepsilon) &\equiv& \int_{R < \|t\| \le (c /(4\bar \rho))^{1/b} n ^{s/(2b)}} \exp\left((n - |\gamma|) \log\left(\frac{1}{n}\sum_{i=1}^n |\phi_{i,n,P}(t)| \right) \right)e^{-(\varepsilon\sqrt{n}\|t\|)^{1/2}} dt \\ \notag
	K_\gamma^3(n,\varepsilon) &\equiv& \int_{(c /(4\bar \rho))^{1/b} n ^{s/(2b)} < \|t\|} \exp\left((n - |\gamma|) \log\left(\frac{1}{n}\sum_{i=1}^n |\phi_{i,n,P}(t)| \right) \right)e^{-(\varepsilon\sqrt{n}\|t\|)^{1/2}} dt.
\end{eqnarray*}
By following the proof of AP, dealing with $K_\gamma^1(n,\varepsilon),K_\gamma^2(n,\varepsilon),K_\gamma^3(n,\varepsilon)$ one by one,  we can show that 
\begin{eqnarray*}
	K_\gamma^1(n,\varepsilon) + K_\gamma^2(n,\varepsilon) + K_\gamma^3(n,\varepsilon)
\end{eqnarray*}
converges to zero faster than any polynomial rate in $n$ uniformly over $P \in \mathscr{P}_n$ even when we take $\varepsilon = n^{-\delta}/2$. Thus we conclude that for any $w \ge 1$, there exist a constant $C$ and a number $n_{0,2} \ge 1$ which depend only on $(s,d,\bar \rho,b,c,R,\delta,w)$ such that for all $n \ge n_{0,2}$,
\begin{eqnarray*}
	I_1 + I_2 + I_3 \le C n^{-w}. 
\end{eqnarray*}
Thus, there exists $n_0 \ge 1$ which depends only on $(s,d,\bar \rho,b,c,R,\delta,w)$ such that for all $n \ge n_0$, the first integral in (\ref{dec}) dominates the second one. Collecting the bounds in (\ref{bound1}),(\ref{bound2}),(\ref{bound3}), and (\ref{bound A}), we find that (taking $\varepsilon = n^{-\delta}/2$)
\begin{eqnarray}
\label{dec2}
|Q_{n,P}(f) - \tilde Q_{n,s,P}(f)| \le c_{12}(s,d) M_s(f) n^{-(s-2)/2} \bar \rho  + \bar \omega_{f_{a_n}}\left(n^{-\delta}; |\tilde Q_{n,s+d,P}'|\right).
\end{eqnarray}
Recall the definition $\tilde Q_{n,s,P}'$ in (\ref{tilde Q'}) and define
\begin{eqnarray*}
	P_j(- \Phi_{0,D_{n,P}};\{\bar \chi_{\nu,P}\}) (A) \equiv \int_A \tilde P_j(-D;\{\bar \chi_{\nu,P}\}) \phi_{0,D_{n,P}}(x)dx, \text{ for any Borel set } A,
\end{eqnarray*}
where $\Phi_{0,D_{n,P}}$ denotes the distribution of $N(0,D_{n,P})$.

As for the last term in (\ref{dec2}), we follow (20.36) and (20.37) of BR, p.213, to find that for absolute constants $C,C'>0$, for any $0 \le j \le s-2$,
\begin{eqnarray*}
	&& \bar \omega_{f_{a_n}}\left( n^{-\delta}; n^{-\frac{j}{2}}|P_j(- \Phi_{0,D_{n,P}};\{\bar \chi_{\nu,P}\})|\right)\\
	&\le& C \rho_{n,s,P}\left( M_s(f)\int(1+\|x\|^{s'}) |\phi_{a_n,D_{n,P}}(x) - \phi(x)|dx + \bar \omega_f(n^{-\delta};\Phi)\right)\\
	&&+ C n^{-\frac{j}{2}}\rho_{n,s,P}\int_{\|x\| > n ^{1/6}}\left(1+\|x\|^{3j+s'}\right)\phi_{a_n,D_{n,P}}(x)dx\\
	&\le& c_{13}(s,d,j) M_s(f) n^{-\frac{s-2}{2}} \bar \Delta_{n,s,P}(1) + C'\rho_{n,s,P}\bar \omega_{f}\left( n^{-\delta};\Phi\right),	
\end{eqnarray*}
where the last inequality uses Lemma 14.6 of BR, p.131. (Here $\phi_{a_n, D_{n,P}}$ denotes the density of $N(a_n,D_{n,P})$.) Also, from (20.39) of BR on p.213, for $s-1 \le j \le s+d-2$,
\begin{eqnarray*}
	&& \bar \omega_{f_{a_n}}\left(n^{-\delta}; n^{-\frac{j}{2}}|P_j(- \Phi_{0,D_{n,P}};\{\bar \chi_{\nu,P}\})|\right)\\
	&\le& c_{14}(s,d,j) n^{-\frac{j}{2}}\frac{1}{n}\sum_{i=1}^n \mathbf{E}_{P}\|Z_{i,n}\|^{j+2}M_s(f)\int\left( 1 + \|x\|^{3j+s}\right) \phi_{a_n,D_{n,P}}(x)dx.
\end{eqnarray*}
From (\ref{bd1}) and (\ref{bd2}), the last term is bounded by
\begin{eqnarray*}
	&& c_{15}(s,d,j) M_s(f) n^{((-s+j+2)-j)/2}\left\{\xi^{-s+j+2} \rho_{n,s,P} + \bar \Delta_{n,s,P}(\xi)\right\}\\
	&\le& c_{16}(s,d,j) M_s(f) n^{-(s-2)/2}\left\{\xi \rho_{n,s,P} + \bar \Delta_{n,s,P}(\xi)\right\},
\end{eqnarray*}
because $0< \xi \le 1$. Hence
\begin{eqnarray*}
	&& \bar \omega_{f_{a_n}}\left( n^{-\delta}; |\tilde Q_{n,s+d,P}'|\right)\\
	&\le& c_{17}(s,d) M_s(f) n^{-(s-2)/2} \left\{\xi \rho_{n,s,P} + \bar \Delta_{n,s,P}(\xi) + \bar \Delta_{n,s,P}(1)\right\} 
	+ C'\rho_{n,s,P}\bar \omega_{f}\left( n^{-\delta};\Phi\right)\\
	&\le& 3 c_{17}(s,d) M_s(f) n^{-(s-2)/2} \bar \rho
	+ C' \bar \rho \bar \omega_{f}\left( n^{-\delta};\Phi\right)
\end{eqnarray*}
Combining this with (\ref{dec2}), we obtain that for all $n \ge n_0$,
\begin{eqnarray*}
		\left| \int f d(Q_{n,P} - \tilde Q_{n,s,P}) \right| \le C n^{-(s-2)/2}\left( 1 + M_s(f)\right) + C \bar \omega_f(n^{-\delta};\Phi),
\end{eqnarray*}
for constants $n_0$ and $C$ that depend only on $s$, $d$, $\bar \rho$, and $(b,c,R)$. The desired bound follows by (\ref{bound}). $\blacksquare$\medskip

We turn to proving Proposition \ref{prop: leb bd}. Define
\begin{eqnarray*}
	\xi(u_i,u_j;t) \equiv \inf_{q \in \mathbb{Z}}\left(t'(u_i - u_j) - 2\pi q \right)^2,
\end{eqnarray*}
and let for $t \in \mathbf{R}^d$,
\begin{eqnarray*}
	A_n(c,r;t) \equiv \left\{ (u_i)_{i=1}^n \in \mathbf{R}^{nd} : \frac{1}{\pi^2 n(n-1)} \sum_{i \ne j}^n\xi(u_i,u_j;t)r^b \le c \right\}.
\end{eqnarray*}

\begin{lemma}
	\label{lemm: inc}
	For any $\tilde t \in \mathbf{R}^d$ such that $\|\tilde t\| >R$,
	\begin{eqnarray*}
		\bigcup_{t \in \mathbf{R}^d:\|t\| > R} A_n(c,\|t\|;t) \subset A_n(c,R;\tilde t).
	\end{eqnarray*}
\end{lemma}

\noindent \textbf{Proof: } Take $t,\tilde t \in \mathbf{R}^d$ such that $\|t\|, \|\tilde t\| >R$, and note that
\begin{eqnarray*}
	t'(u_i - u_j) = \frac{\tilde t' \tilde t t'(u_i - u_j)}{\|\tilde t\|^2} =  \tilde t' (\tilde u_i - \tilde u_j),
\end{eqnarray*}
where $\tilde u_i \equiv \tilde t t' u_i/\|\tilde t\|^2$. As $t'u_i$ runs in $\mathbf{R}$, so does $\tilde t' \tilde u_i$. Hence for any $t_1,t_2 \in \mathbf{R}^d$ such that $\|t_1\|, \|t_2\| > R$, 
\begin{eqnarray*}
	A_n(c,\|t_1\|;t_2) = A_n(c;\|t_1\|,\tilde t).
\end{eqnarray*}
Furthermore, for any $a,a' >0$ such that $a \le a'$,
\begin{eqnarray*}
	A_n(c,a';\tilde t) \subset A_n(c,a;\tilde t).
\end{eqnarray*}
Thus
\begin{eqnarray*}
	\bigcup_{t \in \mathbf{R}^d:\|t\| > R} A_n(c,\|t\|;t) = \bigcup_{t \in \mathbf{R}^d:\|t\| > R} A_n(c,\|t\|;\tilde t) \subset A_n(c,R;\tilde t).
\end{eqnarray*}
$\blacksquare$\medskip

\noindent \textbf{Proof of Proposition \ref{prop: leb bd}:} Let $W$ be a random vector distributed as $n^{-1}\sum_{i=1}^n \delta_{u_i}$. From the proof of Lemma 5.1 of AP, we observe that for all $t \in \mathbf{R}^d$,
\begin{eqnarray*}
	1 - |\mathbf{E}[e^{\text{i}t'W}]|^2 &\ge& \frac{2}{\pi^2 n^2} \sum_{i \ne j}^{n} \xi(u_i,u_j;t) \\
	&=& \frac{2n(n-1)}{\pi^2 n^2} \frac{1}{n(n-1)}\sum_{i \ne j}^{n} \xi(u_i,u_j;t)\\
	&\ge& \frac{2}{\pi^2} \frac{1}{n(n-1)}\sum_{i \ne j}^{n} \xi(u_i,u_j;t).
\end{eqnarray*}
Therefore,
\begin{eqnarray*}
	1 - |\mathbf{E}[e^{\text{i}t'W}]| &\ge& \frac{2}{\pi^2 (1 + |\mathbf{E}[e^{\text{i}t'W}]|)} \frac{1}{n(n-1)} \sum_{i \ne j}^{n} \xi(u_i,u_j;t)\\
	&\ge& \frac{1}{\pi^2 n(n-1)}\sum_{i \ne j}^{n} \xi(u_i,u_j;t).
\end{eqnarray*}
This means that if for all $\|t\| > R$,
\begin{eqnarray*}
	\frac{1}{\pi^2 n(n-1)} \sum_{i \ne j}^{n} \xi(u_i,u_j;t)\|t\|^b \ge c,
\end{eqnarray*}
then $n^{-1}\sum_{i=1}^{n} \delta_{u_i}$ satisfies the weak Cram\'{e}r condition with parameter $(b,c,R)$. Hence for any $c>0$, and any $\|\tilde t\| > R$,
\begin{eqnarray*}
	\mathscr{U}_n^*(b,c,R) \subset \bigcup_{t \in \mathbf{R}^d: \|t\|>R} A_n(c,\|t\|;t) \subset A_n(c,R;\tilde t),
\end{eqnarray*}
by Lemma \ref{lemm: inc}.

Let for $t \in \mathbf{R}^d$,
\begin{eqnarray*}
	c_P(t) \equiv \frac{1}{2\pi^2} \mathbf{E}_P\left[\inf_{q \in \mathbb{Z}}\left( t' (X_{1,n} - X_{2,n}) - 2 \pi q\right)^2\right].
\end{eqnarray*}
We take $\tilde t \in \mathbf{R}^d$ such that $\| \tilde t \| >R$ and for all $n \ge 1$ and $P \in \mathscr{P}_n$,
\begin{eqnarray}
\label{bd26}
c_R \le c_P(\tilde t).
\end{eqnarray}
Such $\tilde t$ exists due to the condition (\ref{bound22}). We find that
\begin{eqnarray}
\label{sub}
\mathscr{U}_n^*(b,c_P(\tilde t)R^b,R) \subset A_n\left(c_P(\tilde t)R^b,R;\tilde t\right).
\end{eqnarray}
We write $P\left\{ X_n \in A_n\left(c_P(\tilde t)R^b,R;\tilde t\right)\right\}$ as
\begin{eqnarray}
\label{prob}
&&P\left\{ \frac{1}{\pi^2 n(n-1)} \sum_{i \ne j}^{n}(\xi_{ij}(\tilde t) - \mathbf{E}[\xi_{ij}(\tilde t)])R^b \le -c_P(\tilde t)R^b \right\}\\ \notag
&=& P\left\{ \frac{1}{n(n-1)} \sum_{i \ne j}^{n}(\xi_{ij}(\tilde t) - \mathbf{E}[\xi_{ij}(\tilde t)]) \le - \pi^2 c_P(\tilde t)  \right\},
\end{eqnarray}
where $\xi_{ij}(t) \equiv \xi(X_{i,n},X_{j,n};t)$. Note that
\begin{eqnarray*}
	\sup_{u_1,u_2 \in \mathbf{R}^d} \xi(u_1,u_2;\tilde t) \le \pi^2.
\end{eqnarray*}
By Proposition 2.3(b) of \cite{Arcones/Gine:93:AP}, the last probability in (\ref{prob}) is less than
\begin{eqnarray*}
	\exp\left( - \frac{ n (\pi^2 c_P(\tilde t))^2}{2 \pi^4}\right) \le \exp\left( - \frac{n  c_P^2(\tilde t)}{2}\right).
\end{eqnarray*}
We conclude that
\begin{eqnarray*}
	P\left\{X_n \in \mathscr{U}_n^*(b,c_P(\tilde t)R^b,R)\right\} \le \exp\left( - \frac{n  c_P^2(\tilde t)}{2}\right).
\end{eqnarray*}
Note that $\mathscr{U}_n^*(b,c,R)$ is increasing in $c$ and $\exp(-nc^2/2)$ is decreasing in $c$. Thus we obtain the desired result from (\ref{bd26}). $\blacksquare$\medskip

\noindent \textbf{Proof of Proposition \ref{prop: suff}:} We show (\ref{bd234}) with $X_{i,n}$ there replaced by $X_{ij,n}$. Suppose that Condition (i) holds. We take $r \in (0,\pi/R)$ and an integer $k$ so that $rk<c_U - c_L < r(k+1)$. Then
\begin{eqnarray*}
	c_P(k;r) \ge \varepsilon^2 \min\left\{(a_{j_2} - a_{j_1} - r k)^2, (a_{j_2} - a_{j_1} - r (k+1))^2\right\},
\end{eqnarray*}
fulfilling condition (\ref{bd234}) for $X_{ij,n}$.

Suppose now that Condition (ii) holds. We take $a'<a''<b''< b'$  and $r \in (0,\pi/R)$ and an integer $k$ such that $a \le a' < a''< b''< b' \le b$ and $0< rk < b'' - a'' < rk + r/2 <  b' -  a' < r(k+1)$. Then, whenever $X_{2,n} \in [a',a'']$ and $X_{1,n} \in [b'',b']$, we have $b'' - a'' \le X_{1,n} - X_{2,n} \le b' - a'$. Hence
\begin{eqnarray*}
	c_P(k;r) &\ge& P\left\{ X_{2,n} \in [a',a'']\right\} P\left\{ X_{1,n} \in [b'',b']\right\} 
	\min\left\{(b'-a' - r(k+1))^2, (b''-a'' - rk)^2\right\}\\
	&\ge& \varepsilon^2 (a''-a')(b'-b'') \min\left\{(b'-a' - r(k+1))^2, (b''-a'' - rk)^2\right\}.
\end{eqnarray*} 
Therefore, again, condition (\ref{bd234}) is fulfilled for $X_{ij,n}$. $\blacksquare$ \medskip

\noindent \textbf{Proof of Theorem \ref{thm: bootstrap}:} Let $\mathscr{U}$ be the set of $u = (u_j)_{j=1}^n \in \mathbf{R}^{2 d n}$ such that the discrete measure $\frac{1}{n}\sum_{j=1}^{n} \delta_{u_j}$ (with $\delta_{u_j}$ denoting Dirac measure at $u_j \in \mathbf{R}^{2d}$) fails to satisfy the weak Cram\'{e}r condition in Definition \ref{def: weak Cramer cond} with parameter $(b,c_R R^b,R)$. We define
\begin{eqnarray*}
	\mathscr{\tilde E}_{n,0} \equiv \left\{(X_{i,n})_{i=1}^{n} \in \mathscr{U}\right\}.
\end{eqnarray*}
By Proposition \ref{prop: leb bd},
\begin{eqnarray*}
	\sup_{ P \in \mathscr{P}_n} P( \mathscr{\tilde E}_{n,0}) \le \exp\left( - \frac{c_R^2 n }{2}\right),
\end{eqnarray*}
which decreases at a rate faster than any polynomial rate in $n$.

Now, we focus on the event $\mathscr{\tilde E}_{n,0}^c\cap \mathscr{E}_n(\bar \rho,c_1)$. Note that by Theorem \ref{thm: Edgeworth} (noting that $\delta > (s-2)/2$) we have
\begin{eqnarray*}
	\left| Q_{n,\mathscr{F}_n}(A) - \tilde Q_{n,\mathscr{F}_n}(A)\right| \le 2C n^{-(s-2)/2} + C \bar \omega_{1\{\cdot \in A\}}\left(n^{-(s-2)/2};\Phi \right),
\end{eqnarray*}
where $C>0$ is the same constant that appears in Theorem \ref{thm: Edgeworth}. Thus we obtain the first statement of the theorem. The second statement follows by Corollary 3.2 of BR, p.24. $\blacksquare$\medskip

\noindent \textbf{Proof of Corollary \ref{cor: mod cont}: } First, we bound
\begin{eqnarray*}
	\left| Q_{n,\mathscr{F}_n}(A^\eta) - Q_{n,\mathscr{F}_n}(A) \right| \le A_{1n} + A_{2n} + A_{3n},
\end{eqnarray*}
where
\begin{eqnarray*}
	A_{1n} &\equiv& \left| Q_{n,\mathscr{F}_n}(A^\eta) - \tilde Q_{n,\mathscr{F}_n}(A^\eta) \right|,\\
	A_{2n} &\equiv& \left| Q_{n,\mathscr{F}_n}(A) - \tilde Q_{n,\mathscr{F}_n}(A) \right|, \text{ and }\\
	A_{3n} &\equiv& \left| \tilde Q_{n,\mathscr{F}_n}(A^\eta) - \tilde Q_{n,\mathscr{F}_n}(A) \right|.
\end{eqnarray*}
By Theorem \ref{thm: bootstrap}, for all $n \ge n_0$,
\begin{eqnarray*}
	&& P\left(\left\{A_{1n} + A_{2n} \le 2(2C + c(s,d))n^{-(s-2)/2} \right\} \cap \mathscr{E}_n(\bar \rho,c_1,c_2) \right) \\
	&& \ge P(\mathscr{E}_n(\bar \rho,c_1,c_2)) - \exp\left( -\frac{c_R^2 n}{2}\right) .
\end{eqnarray*}
On the event $\mathscr{E}_n(\bar \rho,c_1,c_2)$, we have 
\begin{eqnarray*}
	\left| \tilde Q_{n,\mathscr{F}_n}(A^\eta) - \tilde Q_{n,\mathscr{F}_n}(A) \right| \le \left|\tilde Q_{n,\mathscr{F}_n}(A^\eta \setminus A)\right| \le c_1(s,d) c_2 \eta,
\end{eqnarray*}
where $c_1(s,d)>0$ is a constant that depends only on $s,d$. The last inequality follows by Corollary 3.2 of BR, p. 24. Therefore, we obtain the desired result. $\blacksquare$\medskip

\noindent \textbf{Proof of Corollary \ref{cor: bootstrap}: }
Let
\begin{eqnarray*}
	\eta_n \equiv n^{-(s-1)/2} \left(2\sqrt{c_3 d} \log n \right)^s.
\end{eqnarray*}
We follow the arguments on pages 253-254 of \cite{Hall:92:Bootstrap} and write
\begin{eqnarray*}
	\sqrt{n} g_n(\overline X_{n} + n^{-1/2} x ) = g_{n,1}(x) + g_{n,2}(x),
\end{eqnarray*}
where, with $x_k$ denoting the $k$-th entry of $x$, 
\begin{eqnarray*}
	g_{n,1}(x) &\equiv& \sum_{k} \frac{\partial g_n(\overline X_n)}{\partial x_k} x_k + \frac{1}{2}n^{-1/2}\sum_{k_1,k_2} \frac{\partial^2 g_n(\overline X_n)}{\partial x_{k_1} \partial x_{k_2}} x_{k_1} x_{k_2}\\
	&& + ... + \frac{1}{(s-1)!} n^{-(s-2)/2} \sum_{k_1,...,k_{s-1}} \frac{\partial^{s-1} g_n(\overline X_n)}{\partial x_{k_1} ... \partial x_{k_{s-1}}} x_{k_1} \cdots x_{k_{s-1}},
\end{eqnarray*}
and $g_{n,2}(x)$ is such that $|g_{n,2}(x)| \le c'(c_1,s,d) \eta_n$, for some constant $c'(c_1,s,d)>0$ that depends only on $c_1$, $s$ and $d$, whenever $\|x\| \le 2\sqrt{c_3 d} \log n$. Let
\begin{eqnarray*}
	f_{n,t}(x) &\equiv& 1\left\{n^{1/2}g_{n}\left(\overline X_n + n^{-1/2} \hat V_n^{1/2} x \right) \le t \right\},\\
	f_{n,t}^A (x) &\equiv& 1\left\{n^{1/2}g_{n}\left(\overline X_n + n^{-1/2} \hat V_n^{1/2}  x \right) \le t, \|\hat V_n^{1/2} x\| \le 2\sqrt{c_3 d} \log n \right\},\\
	f_{n,t}^B(x) &\equiv& 1\left\{n^{1/2}g_{n}\left(\overline X_n + n^{-1/2} \hat V_n^{1/2}  x \right) \le t, \|\hat V_n^{1/2} x\| > 2 \sqrt{c_3 d} \log n \right\},\text{ and }\\
	f_{n,t}^C(x) &\equiv& 1\left\{g_{n,1}\left(\hat V_n^{1/2} x \right) \le t, \|\hat V_n^{1/2} x\| \le 2\sqrt{c_3 d}  \log n \right\}.
\end{eqnarray*}
Let
\begin{eqnarray*}
	\mathscr{E}_n' \equiv \mathscr{E}_n(\bar \rho,c_1,c_2,c_3),
\end{eqnarray*}
for simplicity. We apply Theorem \ref{thm: bootstrap} to find that for all $n \ge n_0$,
\begin{eqnarray*}
	&& P\left(\left\{\sup_{t \in \mathbf{R}}\left|Q_{n,\mathscr{F}_n}(f_{n,t}) - \tilde Q_{n,\mathscr{F}_n}(f_{n,t})\right| - C \bar \omega_{f_{n,t}}\left(n^{-(s-2)/2};\Phi \right) \le C  n^{-(s-2)/2}\right\} \cap  \mathscr{E}_n' \right)\\
	&& \ge P(\mathscr{E}_n') - \exp\left(-\frac{c_R^2n}{2} \right),
\end{eqnarray*}
where $C$ and $n_0$ are constants that appear in Theorem \ref{thm: bootstrap}. Observe that
\begin{eqnarray*}
	\bar \omega_{f_{n,t}}\left(n^{-(s-2)/2};\Phi \right) &\le& \bar \omega_{f_{n,t}^A}\left(n^{-(s-2)/2};\Phi \right) + \bar \omega_{f_{n,t}^B}\left(n^{-(s-2)/2};\Phi \right)\\
	&\le& \bar \omega_{f_{n,t}^C}\left(n^{-(s-2)/2} + c'(c_1,s,d) \eta_n ;\Phi \right) + \bar \omega_{f_{n,t}^B}\left(n^{-(s-2)/2};\Phi \right).
\end{eqnarray*}
By Lemma 5.3 of \cite{Hall:92:Bootstrap}, we have
\begin{eqnarray*}
	\bar \omega_{f_{n,t}^C}\left(n^{-(s-2)/2} + c'(c_1,s,d) \eta_n ;\Phi \right) &\le& c''(c_1,s,d) \left(n^{-(s-2)/2} + \eta_n \right) \\
	&\le& 2 c''(c_1,s,d) n^{-(s-2)/2},
\end{eqnarray*}
for all $n \ge n_1$, where $c''(c_1,s,d)$ depends only on $c_1,s,d$ and $n_1$ depends only on $s$, $d$ and $c_3$.

Note that on the event $\mathscr{E}_n'$,
\begin{eqnarray}
\label{ineq2}
	\|\hat V_n^{1/2} x\| \le \|\hat V_n^{1/2}\|\|x\| \le \sqrt{\text{tr}(\hat V_n)} \|x\|
	\le \sqrt{c_3 d} \|x\|. 
\end{eqnarray}
Also, observe that for any real functions $a(x)$ and $b(x)$ of $x \in \mathbf{R}^d$,
\begin{eqnarray*}
	&& \left|1\{a(x) \le t, b(x) > c\} - 1\{a(x+y) \le t, b(x+y) > c\}\right|\\
	&\le& \left|1\{a(x) \le t, b(x) > c\} - 1\{a(x+y) \le t, b(x) > c\}\right|\\
	&& + \left|1\{a(x+y) \le t, b(x) > c\} - 1\{a(x+y) \le t, b(x+y) > c\}\right|\\
	&\le& 1\{b(x) > c\} + \left|1\{b(x) > c\} - 1\{b(x+y) > c\}\right|.
\end{eqnarray*}
Hence, $\bar \omega_{f_{n,t}^B}(n^{-(s-2)/2};\Phi)$ is bounded by
\begin{eqnarray*}
    && \int  1\left\{2 \sqrt{c_3 d} \log n < \|\hat V_n^{1/2} x\| \right\} d\Phi(x)\\
    && + \int  1\left\{2 \sqrt{c_3 d} \log n - n^{-(s-2)/2} < \|\hat V_n^{1/2} x\| < 2 \sqrt{c_3 d} \log n + n^{-(s-2)/2}\right\} d\Phi(x) \\
   &\le& 2\int  1\left\{\sqrt{c_3 d} \log n  < \|\hat V_n^{1/2} x\| \right\} d\Phi(x)
   \le 2 \int_{x: \|x\| \ge \log n }d\Phi,
\end{eqnarray*}
for all $n \ge n_2$, where $n_2 \ge 1$ is such that $\sqrt{c_3d}\log n_2 > n_2^{-(s-2)/2}$ and $n_2$ depends only on $s$, $c_3$ and $d$. The last inequality follows from (\ref{ineq2}). By Markov's inequality, the last integral is bounded by
\begin{eqnarray}
\label{eq}
	\exp\left(-\frac{(\log n)^2}{4}\right)\int \exp\left(\frac{\| x\|^2}{4}\right)d \Phi(x)
	= n^{-\frac{1}{4}\log n} \left(\frac{1}{2}\right)^{-d/2},
\end{eqnarray}
as the integral on the left hand side is the moment generating function of $\chi_d^2$ at $1/4$. We can take $n_3 \ge 1$ that depends only on $s,d$ such that the last term in (\ref{eq}) is bounded by $n^{-(s-2)/2}$. Taking $n_0'$ to be the maximum of $n_0,n_1,n_2,n_3$, we obtain the corollary. $\blacksquare$\medskip

\bibliographystyle{econometrica}
\bibliography{uniformEdgeworth}
\end{document}